\documentclass[12pt,a4paper]{article}
\usepackage[cp1251]{inputenc} 
\usepackage{amsfonts}

\newcommand{\di}{\displaystyle}
\newcommand{\B}{$\hfill\Box$}
\newcommand{\al}{\alpha}
\newcommand{\be}{\beta}
\newcommand{\ga}{\gamma}
\newcommand{\de}{\delta}
\newcommand{\la}{\lambda}
\newcommand{\om}{\omega}
\newcommand{\ee}{\varepsilon}

\newcommand{\iy}{\infty}

\begin{document}

\begin{center}
{\large\bf
Recovering Differential Operators with Nonseparated Boundary Conditions
in the Central Symmetric Case}\\[0.1cm]
{\bf V.\,Yurko} \\[0.1cm]
\end{center}

\thispagestyle{empty}

{\bf Abstract.} Inverse spectral problems for Sturm-Liouville operators on
a finite interval with non-separated boundary conditions are studied in
the central symmetric case, when the potential is symmetric with respect
to the middle of the interval. We discuss statements of the problems,
provide algorithms for their solutions along with necessary and sufficient
conditions for the solvability of the inverse problems considered.



\noindent {\bf 1. Introduction. }\\

We study inverse spectral problems for the Sturm-Liouville operator
$$
\ell y:=y''+q(x)y,\; x\in (0,\pi),
$$
on the finite interval $(0,\pi)$ with non-separated boundary conditions.
Inverse problems consist in recovering coefficients of differential
operators from their spectral characteristics. Such problems often appear
in mathematics, mechanics, physics, geophysics, electronics and other
branches of natural sciences and engineering. Inverse problems also play
an important role in solving nonlinear evolution equations in mathematical
physics. Inverse problems for differential operators with separated boundary
conditions have been studied fairly completely by many authors (see the
monographs [1-5] and the references therein). Inverse problems for
Sturm-Liouville operators with non-separated boundary conditions, which
are more difficult for the investigation, were treated in [6-17] and
other works. In particular, the periodic boundary value problem was
considered in [6, 7, 9, 14]. Stankevich [6] suggested a statement of
the inverse problem and proved the corresponding uniqueness theorem.
Marchenko and Ostrovskii [7] gave the characterization of the spectrum
for the periodic boundary value problem in terms of a special conformal
mapping. Conditions considered in [7] are difficult to verify. Another
method, used in [9], allowed to obtain necessary and sufficient conditions
for the solvability of the inverse problem for the periodic case that
are easier to verify. Similar results were obtained in [9] for another
type of boundary conditions, namely
$$
y'(0)-ay(0)+by(\pi)=y'(\pi)+dy(\pi)-by(0)=0.
$$
Later analogous results were established in [12-13].

In this paper we study the case when the potential $q$ is symmetric with
respect to the middle of the interval, i.e. $q(x)=q(\pi-x)$ a.e. on
$(0,\pi).$ The symmetric case requires nontrivial modifications in the
method and allows us to specify less spectral information than in the
general case. Some results for the symmetric case were obtained in [10]
and [17]. In the present paper for the symmetric case we construct the
solution of the inverse spectral problem and give the characterization
of the spectrum for various non-separated boundary conditions. For
convenience of readers in Section 2 we describe briefly the known
results for the general (non-symmetric) case.\\

{\bf 2. Periodic boundary value problem.}\\

Consider the differential equation
$$
-y''+q(x)y=\la y,\quad x\in (0,\pi),                                 \eqno(1)
$$
where $\la$ is the spectral parameter, and $q(x)\in L_2(0,T)$ is a
real-valued function. The function $q(x)$ is called the potential.
Let $C(x,\la), S(x,\la)$ and $\psi(x,\la)$ be solutions of Eq. (1)
with the initial conditions $C(0,\la)=S'(0,\la)=-\psi'(\pi,\la)=1,$
$C'(0,\la)=S(0,\la)=\psi(\pi,\la)=0.$ For each fixed $x,$ the
functions $C^{(\nu)}(x,\la), S^{(\nu)}(x,\la)$ and $\psi^{(\nu)}(x,\la),
\nu=0,1,$ are entire in $\la$ of order $1/2.$
Moreover,
$$
\langle C(x,\la), S(x,\la)\rangle\equiv 1,                           \eqno(2)
$$
where $\langle y,z\rangle :=yz'-y'z$ is the Wronskian of $y$ and $z.$
Denote
$$
\Delta(\la)=(C(\pi,\la)+S'(\pi,\la))/2,\;
\delta(\la)=(C(\pi,\la)-S'(\pi,\la))/2,\; p(\la)=1-\Delta(\la).
$$
Zeros $\Lambda=\{\la_n\}_{n\ge 0}$ of the entire function $p(\la)$
coincide with the eigenvalues of the boundary value problem (BVP)
$L=L(q)$ for Eq. (1) with periodic boundary conditions
$$
y(0)-y(\pi)=y'(0)-y'(\pi)=0.
$$
The function $p(\la)$ is called the characteristic function for $L.$
For convenience of readers we describe briefly the well-known results
related to the BVP $L$ (see [6, 7, 9] for details).

1) All eigenvalues $\la_n$ are real, and
$$
\la_0<\la_1\le\la_2<\la_3\le\la_4<\ldots,                           \eqno(3)
$$
$$
\la_{2n}=(2n)^2+\al+\kappa_{2n},\;
\la_{2n-1}=(2n)^2+\al+\kappa_{2n-1},\quad \{\kappa_n\}\in l_2,      \eqno(4)
$$
where $\al=\frac{1}{\pi}\int_0^\pi q(t)\,dt.$ Here and everywhere
below one and the same symbol $\{\kappa_n\}$ denotes various
sequences from $l_2$. The specification of $\Lambda$ uniquely
determines the characteristic function $p(\la)$ by the formula
$$
p(\la)=\frac{\pi^2}{2}(\la-\la_0)
\prod_{n=1}^\iy \frac{\la_{2n}-\la}{(2n)^2}
\prod_{n=1}^\iy \frac{\la_{2n-1}-\la}{(2n)^2}.                      \eqno(5)
$$
Moreover,
$$
\max_{\la\in[\la_{2n},\la_{2n+1}]} p(\la)\ge 2,\quad n\ge 0.       \eqno(6)
$$

2) Let $\Lambda^{+}=\{\la^{+}_n\}_{n\ge 1}$ be zeros of the entire
function $p^{+}(\la):=p(\la)-2.$ Then $\{\la^{+}_n\}_{n\ge 0}$
are real and
$$
\la_0<\la^{+}_1\le\la^{+}_2<\la_1\le\la_2<
\la^{+}_3\le\la^{+}_4<\la_3\le\la_4\ldots,                         \eqno(7)
$$
$$
\la^{+}_{2n}=(2n-1)^2+\al+\kappa_{2n},\;
\la^{+}_{2n-1}=(2n-1)^2+\al+\kappa_{2n-1},\quad\{\kappa_n\}\in l_2. \eqno(8)
$$
Denote $a_{2n}=[\la_{2n-1},\la_{2n}],\;a_{2n-1}=
[\la^{+}_{2n-1},\la^{+}_{2n}],\; n\ge 1.$ Segments $a_n$ are
called the gaps.

\smallskip
3) Denote $d(\la):=\langle \psi(x,\la), S(x,\la)\rangle
=S(\pi,\la)=\psi(0,\la).$ Then zeros $\ga=\{\ga_n\}_{n\ge 1}$ of the
entire function $d(\la)$ coincide with the eigenvalues of the BVP
$L_0=L_0(q)$ for Eq. (1) with Dirichlet boundary conditions
$y(0)=y(\pi)=0.$ The numbers $\ga_n$ are real, $\ga_n\in a_n$, and
$$
\ga_1<\ga_2<\ga_3<\ldots;\quad
\ga_n=n^2+\al+\kappa_n,\; \{\kappa_n\}\in l_2.                     \eqno(9)
$$
The specification of $\ga$ uniquely determines the characteristic
function $d(\la)$ of $L_0$ by the formula
$$
d(\la)=\pi\prod_{n=1}^\iy \frac{\ga_n-\la}{n^2}\,.                \eqno(10)
$$
The numbers $\al_n:=\int_0^\pi S^2(x,\ga_n)\,dx$ are called the
weight numbers, and numbers $\{\ga_n,\al_n\}_{n\ge 1}$ are called
the spectral data for the BVP $L_0$. One has
$$
\al_n=\dot d(\ga_n)S'(\pi,\ga_n),
\quad \dot d(\la):=\frac{d}{d\la}d(\la),                          \eqno(11)
$$
$$
\al_n>0;\qquad \al_n=\frac{\pi}{2n^2}
\Big(1+\frac{\kappa_n}{n}\Big),\quad \{\kappa_n\}\in l_2,         \eqno(12)
$$
$$
\dot d(\ga_n)=\frac{(-1)^n\pi}{2n^2}
\Big(1+\frac{\kappa_n}{n}\Big),\quad \{\kappa_n\}\in l_2,
\qquad \mbox{sign}\,\dot d(\ga_n)=(-1)^n.                        \eqno(13)
$$
The functions $S(x,\ga_n)$ and $\psi(x,\ga_n)$ are eigenfunctions
for $L_0$, and
$$
\psi(x,\ga_n)=\be_n S(x,\ga_n),\quad \be_n\ne 0.                 \eqno(14)
$$

{\bf Lemma 1. }{\it The following relation holds}
$$
\al_n\be_n=-\dot d(\ga_n).                                       \eqno(15)
$$

{\it Proof. } Since
$$
-\psi''(x,\lambda)+q(x)\psi(x,\lambda)=\lambda\psi(x,\lambda), \;
-S''(x,\ga_n)+q(x)S(x,\ga_n)=
\ga_n S(x,\ga_n),
$$
we get
$$
\di\frac{d}{dx}\langle\psi(x,\lambda),S(x,\ga_n)\rangle=
(\lambda-\ga_n)\psi(x,\lambda)S(x,\ga_n),
$$
and hence,
$$
(\lambda-\ga_n)\di\int_0^\pi \psi(x,\lambda)S(x,\ga_n)\,dx=
\langle\psi(x,\lambda),S(x,\ga_n)\rangle\Big|_0^\pi=-d(\lambda).
$$
For $\lambda \to \lambda_n,$ this yields
$$
\di\int_0^\pi\psi(x,\ga_n)S(x,\ga_n)\,dx=
-\dot d(\ga_n).
$$
Using (14) we arrive at (15).
\B

\smallskip
The inverse problem for the BVP $L_0$ is formulated as follows.

\smallskip
{\bf Inverse problem 1. } Given the spectral data
$\{\ga_n,\al_n\}_{n\ge 1}$, construct the potential $q.$

\smallskip
This inverse problem is related to the case of separated boundary
conditions. It is known that the specification of the spectral
data $\{\ga_n,\al_n\}_{n\ge 1}$ uniquely determines the potential
$q.$ The global solution of Inverse problem 1 can be constructed
by the transformation operator method or by the method of spectral
mappings (see [1-5] for details). In particular, these methods
allow one to describe necessary and sufficient conditions for the
solvability of Inverse problem 1 which are presented in the next
theorem.

\smallskip
{\bf Theorem 1. }{\it For real numbers $\{\ga_n,\al_n\}_{n\ge 1}$ to
be the spectral data for a certain BVP $L_0$ with a real potential
$q(x)\in L_2(0,\pi),$ it is necessary and sufficient that (9)
and (12) hold.}

\smallskip
Let us now return to the periodic BVP $L.$
It follows from (2) that
$$
\Delta^2(\la)-\de^2(\la)-d(\la)d_1(\la)\equiv 1,                 \eqno(16)
$$
where $d_1(\la):=C'(\pi,\la).$ In particular, (16) yields
$$
\de^2(\ga_n)=\Delta^2(\ga_n)-1.                                  \eqno(17)
$$
Denote $\Omega=\{\om_n\}_{n\ge 1},\;\om_n=\mbox{sign}\,\de(\ga_n).$
The sequence $\Omega$ is called the $\Omega$- sequence for $q.$
In view of (17) one has
$$
\de(\ga_n)=\om_n(\Delta^2(\ga_n)-1)^{1/2},                       \eqno(18)
$$
Since $S'(\pi,\ga_n)=\Delta(\ga_n)-\de(\ga_n),$
it follows from (11) and (18) that
$$
\al_n=
\dot d(\ga_n)(\Delta(\ga_n)-\om_n(\Delta^2(\ga_n)-1)^{1/2}).     \eqno(19)
$$

The inverse problem for the periodic case is formulated as
follows [6].

\smallskip
{\bf Inverse problem 2. } Given $\Lambda, \ga$ and $\Omega,$
construct the potential $q.$

\smallskip
This inverse problem was studied in [6, 7, 9, 14] and other works.
It was proved in [6] that the specification of $\Lambda, \ga$
and $\Omega$ uniquely determines the potential $q.$ In order to
construct $q$ one can calculate the functions $p(\la)$ and $d(\la)$
according to (5) and (10), and construct $\{\al_n\}_{n\ge 1}$ via
(19), where $\Delta(\la)=1-p(\la).$ Then using data
$\{\ga_n,\al_n\}_{n\ge 1}$, we can construct the potential $q$
by solving Inverse problem 1.

\smallskip
{\bf Lemma 2. }{\it Fix $n\ge 1.$ Relation $\de(\ga_n)=0$ holds iff
$\ga_n$ lies at one of the endpoints of the gap $a_n$.}

\smallskip
Indeed, in view of (17), $\de(\ga_n)=0,$ iff $\Delta(\ga_n)=\pm 1,$
i.e. $\ga_n$ lies at one of the endpoints of the gap $a_n$.

\smallskip
Denote by $J$ the set of sequences $\Omega=\{\om_n\}_{n\ge 1}$
such that $\om_n=0$ if $\ga_n$ lies at one of the endpoints of
the gap $a_n$, and $\om_n=\pm 1,$ otherwise. Clearly, if
$\Omega$ is the $\Omega$- sequence for $L,$ then $\Omega\in J.$
The next theorem [9] establishes necessary and sufficient
conditions for the solvability of Inverse problem 2.

\smallskip
{\bf Theorem 2 [9]. }{\it Let real numbers
$\Lambda=\{\la_n\}_{n\ge 0}$ satisfying (3)-(4) be given.
The sequence $\Lambda$ is the spectrum for a certain BVP $L$ with
a real potential $q(x)\in L_2(0,\pi),$ iff relation (6) holds, where
$p(\la)$ is constructed via (5). Moreover, if additionally we have
a sequence $\ga=\{\ga_n\}_{n\ge 1}$, $\ga_n\in a_n$, satisfying (9),
where $\Lambda^{+}=\{\la^{+}_n\}_{n\ge 1}$ are zeros of
$p^{+}(\la)=p(\la)-2,$ and a sequence $\Omega=\{\om_n\}_{n\ge 1}\in J,$
then there exists a unique real function $q(x)\in L_2(0,\pi)$ such
that $\Lambda$ and $\ga$ are the spectra of $L$ and $L_0$,
respectively, and $\Omega$ is the $\Omega$- sequence for $L.$}

\smallskip
The next theorem [9] shows that one of the endpoints of each gap can
be chosen arbitrary taking only asymptotics into account.

\smallskip
{\bf Theorem 3 [9]. }{\it Let real numbers $\theta_n$ of the form
$\theta_n=n^2+\al+\kappa_n,\; \{\kappa_n\}\in l_2$,
$\theta_n<\theta_{n+1},$ be given. Then there exists a real function
$q(x)\in L_2(0,\pi)$ (not unique!) such that for this potential the number
$\theta_n$ lies at one of the endpoints of the gap $a_n$ for all $n\ge 1.$}

\bigskip
{\bf 3. Central symmetric case.}\\

In this section we consider the case when the potential $q$ is symmetric
with respect to the middle of the interval, i.e. with respect to the
replacement $x\to \pi-x.$ We will say that $q(x)\in L'_2(0,\pi)$ if
$q(x)\in L_2(0,\pi)$ and $q(x)=q(\pi-x)$ a.e. on $(0,\pi).$

\smallskip
{\bf Theorem 4. }{\it $q(x)\in L'_2(0,\pi)$ iff $\be_n=(-1)^{n-1},\;n\ge 1.$}

\smallskip
{\it Proof. } 1) Let $q(x)\in L'_2(0,\pi).$ Then
$\psi(x,\la)\equiv S(\pi-x,\la).$ Using (14) we calculate
$$
\psi(x,\ga_n)=\be_n S(x,\ga_n)=\be_n\psi(\pi-x,\ga_n)=
\be_n^2 S(\pi-x,\ga_n)=\be_n^2\psi(x,\ga_n).
$$
Hence, $\be_n^2=1.$ On the other hand, it follows from (14) that
$\be_n S'(\pi,\ga_n)=-1.$ Using Sturm's oscillation theorem we
conclude that $\be_n=(-1)^{n-1},\;n\ge 1.$

2) Let $\be_n=(-1)^{n-1},\;n\ge 1.$ Denote $\tilde q(x):=q(\pi-x).$
We agree that here and below, if a certain symbol $\theta$ denotes
an object related to $q,$ then $\tilde\theta$ will denote the
analogous object related to $\tilde q.$

Obviously, $\tilde\psi(x,\la)\equiv S(\pi-x,\la),$
$\tilde S(x,\la)\equiv\psi(\pi-x,\la),$ and consequently,
$d(\la)\equiv\tilde d(\la)$ and $\ga_n=\tilde\ga_n$, $n\ge 1.$
Since $\be_n=(-1)^{n-1},$ it follows from (14) that
$\psi(x,\ga_n)=(-1)^{n-1} S(x,\ga_n).$ Moreover, according to
(14), $\tilde\psi(x,\ga_n)=\tilde\be_n \tilde S(x,\ga_n),$
hence $S(\pi-x,\ga_n)=\tilde\be_n\psi(\pi-x,\ga_n),$ i.e.
$\tilde\be_n=(\be_n)^{-1}=(-1)^{n-1}.$
Thus, $\be_n=\tilde\be_n$ for all $n\ge 1.$ Taking (15) into
account we conclude that $\al_n=\tilde\al_n$ for all $n\ge 1.$
Since the specification of the spectral data $\{\ga_n,\al_n\}_{n\ge 1}$
uniquely determines the potential, we obtain that $q(x)=\tilde q(x)$
a.e. on $(0,\pi),$ i.e. $q(x)\in L'_2(0,\pi).$
\B

\smallskip
Let us consider the inverse problem for the BVP $L_0$. In the central
symmetric case $q(x)\in L'_2(0,\pi)$ we do not need to specify the
weight numbers $\{\al_n\}_{n\ge 1}$; it is sufficient to specify
only the spectrum $\ga.$

\smallskip
{\bf Inverse problem 3. } Given the spectrum $\ga=\{\ga_n\}_{n\ge 1}$,
construct the potential $q.$

\smallskip
It is known [1-5] that for the central symmetric case the
specification of the spectrum $\ga=\{\ga_n\}_{n\ge 1}$ of the BVP
$L_0$ uniquely determines the potential $q.$ In order to construct
$q,$ one can calculate $d(\la)$ via (10) and the weight numbers
$\al_n=(-1)^n \dot d(\ga_n),$ and then find $q$ by solving
Inverse problem 1. Moreover, the characterization of
the spectrum of $L_0$ is given by the following assertion.

\smallskip
{\bf Theorem 5. }{\it For real numbers $\{\ga_n\}_{n\ge 1}$
to be the spectrum of a BVP $L_0$ with a real potential
$q(x)\in L'_2(0,\pi),$ it is necessary and sufficient that (9) holds.}

\smallskip
{\it Proof. } The necessity is obvious. We will prove the sufficiency.
Let real numbers $\{\ga_n\}_{n\ge 1}$ satisfying (9) be given.
We construct $d(\la)$ via (10) and the numbers $\{\al_n\}_{n\ge 1}$
by $\al_n=(-1)^n \dot d(\ga_n).$
Our plan is to use Theorem 1. For this purpose we should obtain the
asymptotics for the numbers $\al_n.$ This seems to be difficult because
the function $d(\la)$ is by construction the infinite product. But
fortunately, for calculating the asymptotics of $\al_n$ one can also
use Theorem 1, as an auxiliary assertion. Indeed, by virtue of Theorem 1
there exists a potential $\tilde q(x)\in L_2(0,\pi)$ (not unique) such
that $\ga=\{\ga_n\}_{n\ge 1}$ is the spectrum of $\tilde L_0:=L_0(\tilde q)$
with this potential. Then $d(\la)$ is the characteristic function of
$\tilde L_0$, and consequently, (13) holds. Therefore, (12) is valid.
Then, by Theorem 1 there exists a unique potential $q(x)\in L_2(0,\pi)$
such that $\{\ga_n,\,\al_n\}_{n\ge 1}$ are the spectral data of $L_0(q).$
Since $\be_n=(-1)^{n-1},\; n\ge 1,$ it follows from Theorem 4 that
$q(x)\in L'_2(0,\pi).$
\B

\smallskip
{\bf Theorem 6 [9]. }{\it $q(x)\in L'_2(0,\pi)$ iff $\ga_n$ lies at
one of the endpoints of the gap $a_n$ for all $n\ge 1.$}

\smallskip
{\it Proof. } 1) Let $q(x)=q(\pi-x)$ a.e. on $(0,\pi).$ Using
Lemma 4 from [8] we get $C(\pi,\la)\equiv S'(\pi,\la),$ i.e.
$\de(\la)\equiv 0.$ By Lemma 2 we conclude that $\ga_n$ lies at one
of the endpoints of the gap $a_n$ for all $n\ge 1.$

2) Let $\ga_n$ lie at one of the endpoints of the gap $a_n$ for all
$n\ge 1.$ By Lemma 1 one has $\de(\ga_n)=0$ for all $n\ge 1.$ Then
the function $F(\la):=\de(\la)/d(\la)$ is entire in $\la,$ and it
vanishes at infinity. This means that $F(\la)\equiv 0,$ and
consequently, $C(\pi,\la)\equiv S'(\pi,\la).$ Using Lemma 4
from [8] we get $q(x)=q(\pi-x)$ a.e. on $(0,\pi).$
\B

\smallskip
We will write $a_n\in I_0$, if the length of the gap $a_n$ is equal
to zero, and $a_n\in I_1$, otherwise.

\smallskip
Let us now consider the inverse problem for the periodic BVP $L.$
In the general case in Inverse problem 2 we have to specify
$\Lambda,\,\ga$ and $\Omega.$ In the central symmetric case we do not
need $\ga.$ On the other hand, the sequence $\Omega=\{\om_n\}_{n\ge 1}$
does not bring any information because in the central symmetric case
$\om_n=0$ for all $n\ge 1.$ Unfortunately, in contrast to the separated
boundary conditions, for the periodic case the specification of the
spectrum $\Lambda$ does not uniquely determine the potential $q,$
and we need additional information. For this purpose we
introduce the sequence $E=\{\ee_n\}_{n\ge 1}$, where $\ee_n=0$ if
$a_n\in I_0$, $\ee_n=1,$ if $a_n\in I_1$ and $\ga_n$ lies at the right
endpoint of $a_n$, $\ee_n=-1,$ if $a_n\in I_1$ and $\ga_n$ lies at the
left endpoint of $a_n$. The sequence $E=\{\ee_n\}_{n\ge 1}$ is called
the $E$- sequence for the potential $q(x)\in L'_2(0,\pi).$
The inverse problem for the periodic BVP $L$ in the central symmetric
case is formulated as follows.

\smallskip
{\bf Inverse problem 4. } Given $\Lambda$ and $E,$ construct $q.$

\smallskip
{\bf Theorem 7 [9]. }{\it Let $q(x)\in L'_2(0,\pi).$ Then the
specification of $\Lambda$ and $E$ uniquely determines the potential $q.$
The solution of Inverse problem 1 can be found by the following algorithm.}

{\bf Algorithm 1. } {\it Given $\Lambda$ and $E.$\\
1) Construct $p(\la)$ by (5).\\
2) Calculate the functions $\Delta(\la)=1-p(\la)$ and $p^{+}(\la)=p(\la)-2.$\\
3) Find zeros $\Lambda^{+}=\{\la^{+}_n\}_{n\ge 1}$ of $p^{+}(\la).$\\
4) Construct $\ga=\{\ga_n\}_{n\ge 1}$ as follows:
$\ga_n$ lies at the right endpoint of $a_n$ if $\ee_n=1;$
$\ga_n$ lies at the left endpoint of $a_n$ if $\ee_n=-1,$ and
$\ga_n=a_n$ if $\ee=0.$\\
5) Using $\{\ga_n\}$ construct the potential $q(x)\in L'_2(0,\pi)$
by solving Inverse problem 3.}

\smallskip
Denote by $J_1$ the set of sequences $E=\{\ee_n\}_{n\ge 1}$
such that $\ee_n=0$ if $a_n\in I_0$, and $\ee_n=\pm 1$ if $a_n\in I_1$.
Clearly, if $E$ is the $E$- sequence for $q,$ then $E\in J_1.$
The next theorem [9] establishes necessary and sufficient
conditions for the solvability of Inverse problem 4.

\smallskip
{\bf Theorem 8 [9]. }{\it Let real numbers
$\Lambda=\{\la_n\}_{n\ge 0}$ satisfying (3)-(4) be given.
The sequence $\Lambda$ is the spectrum for a certain BVP $L$
with a real potential $q(x)\in L'_2(0,\pi),$ iff relation (6) holds,
where $p(\la)$ is constructed via (5). Moreover, if additionally
we have a sequence $E=\{\ee_n\}_{n\ge 1}\in J_1,$ then there exists
a unique real function $q(x)\in L'_2(0,\pi)$ such that $\Lambda$
is the spectrum of $L,$ and $E$ is the $E$- sequence for $q.$}

\smallskip
{\it Proof. } The necessity is obvious. We will prove the sufficiency.
Let real numbers $\Lambda=\{\la_n\}_{n\ge 0}$ satisfying (3)-(4) be given.
We construct the function $p(\la)$ by (5), and calculate the functions
$\Delta(\la)=1-p(\la)$ and $p^{+}(\la)=p(\la)-2.$ Let (6) holds.
Then there exist zeros $\Lambda^{+}=\{\la^{+}_n\}_{n\ge 1}$ of the
function $p^{+}(\la)$, and (7) holds. Using (5) by similar arguments
as in the proof of Theorem 5 (see also [NOVA, p.45]) one gets
$$
p(\la)=1-\cos\rho\pi-
a\frac{\sin\rho\pi}{\rho}-\frac{\kappa(\rho)}{\rho}\,,            \eqno(20)
$$
where $\kappa(\rho)\in L_2(-\iy,\iy)$ for real $\rho.$
Since $p^{+}(\la)=p(\la)-2,$ it follows from (20) that (8) is valid.
Let a sequence $E=\{\ee_n\}_{n\ge 1}\in J_1$ be given. We introduce
real numbers $\ga=\{\ga_n\}_{n\ge 1}$ as follows:
$\ga_n$ lies at the right endpoint of $a_n$ if $\ee_n=1;$
$\ga_n$ lies at the left endpoint of $a_n$ if $\ee_n=-1;$
$\ga_n=a_n$ if $\ee_n=0.$ Clearly, (9) is valid.
We construct the function $d(\la)$ by (10), and the sequence
$\{\al_n\}_{n\ge 1}$ via
$$
\al_n=\dot d(\ga_n)\Delta(\ga_n),\quad n\ge 1.                     \eqno(21)
$$
Since $\Delta(\la)=1-p(\la),$ it follows from (20) that
$$
\Delta(\la)=
\cos\rho\pi+a\frac{\sin\rho\pi}{\rho}+\frac{\kappa(\rho)}{\rho}\,.
$$
Together with (9) this yields
$$
\Delta(\ga_n)=(-1)^n
\Big(1+\frac{\kappa_n}{n}\Big),\quad \{\kappa_n\}\in l_2.         \eqno(22)
$$
Moreover, (13) is valid. It follows from (13), (21) and (22) that
(12) holds. It is easy to check that
$$
\mbox{sign}\,\dot d(\ga_n)=(-1)^n,
\quad \mbox{sign}\,\Delta(\ga_n)=(-1)^n.                         \eqno(23)
$$
In view of (21) and (23) we conclude that $\al_n>0,\; n\ge 1.$
By Theorem 1 we infer that there exists a unique real potential
$q(x)\in L_2(0,\pi)$ such that $\{\ga_n,\al_n\}_{n\ge 1}$
are the spectral data for the BVP $L_0$ for this potential.
We construct solutions $C(x,\la), S(x,\la)$ for Eq.(1) with
this potential. Denote
$$
\tilde\Delta(\la)=(C(\pi,\la)+S'(\pi,\la))/2,\;\tilde p(\la)=
1-\tilde\Delta(\la),\; \tilde p^{+}(\la)=\tilde p(\la)-2.
$$
Using (10) and (21) we get
$$
\Delta(\ga_n)=\tilde\Delta(\ga_n),\quad n\ge 1.
$$
Then the function $F_0(\la):=(\Delta(\la)-\tilde\Delta(\la))/d(\la)$
is entire in $\la,$ and it vanishes at infinity. This yields
$F_0(\la)\equiv 0,$ i.e. $\Delta(\la)\equiv\tilde\Delta(\la),$
and consequently, $p(\la)\equiv \tilde p(\la),$
$p^{+}(\la)\equiv \tilde p^{+}(\la).$ In particular, this means that
the sequence $\Lambda=\{\la_n\}_{n\ge 0}$ coincides with the spectrum
of the BVP $L$ for the potential $q.$ Since $\ga_n$ lies at one of
the endpoints of the gap $a_n$ for all $n\ge 1,$ it follows from
Theorem 6 that $q(x)\in L'_2(0,\pi).$ Now it is clear that $E$ is
the $E$- sequence for $q.$
\B

\medskip
Similar results are valid for other non-separated boundary conditions.
For convenience of readers and for completeness of the presentation,
we formulate here briefly the main results from [10] related to the
boundary conditions
$$
y'(0)-ay(0)+by(\pi)=y'(\pi)+ay(\pi)-by(0)=0.                      \eqno(24)
$$
We consider the BVP $B$ for the differential equation
$$
-y''+q(x)y=\la y,\; x\in (0,\pi),\; q(x)\in L'_2(0,\pi),          \eqno(25)
$$
with the non-separated boundary conditions (24),
where $a$ and $b$ are real numbers, $b\ne 0.$
Let $\theta(x,\la)$ be the solution of Eq. (25) under the initial
conditions $\theta(0,\la)=1,\; \theta'(0,\la)=a.$ Eigenvalues
$\mu=\{\mu_n\}_{n\ge 0}$ of the BVP $B$ coincide with zeros of
the entire function
$$
r(\la)=-\theta'(\pi,\la)-a\theta(\pi,\la)+b^2 S(\pi,\la)+2b.
$$
The eigenvalues $\mu_n$ are real, and
$$
\mu_n<\mu_{n+2},\quad \mu_n=
n^2+\pi^{-1}(h+(-1)^{n+1}4b)+\kappa_n,\; \{\kappa_n\}\in l_2,   \eqno(26)
$$
where $h=4a+\int_0^\pi q(t)\,dt.$ The specification of the
spectrum $\mu$ uniquely determines the characteristic function
$r(\la)$ via
$$
r(\la)=\pi(\la-\mu_0)\prod_{n=1}^\iy \frac{\mu_n-\la}{n^2}\,.   \eqno(27)
$$
Moreover,
$$
\max_{\la\in Q_n} |r(\la)|\ge |4b|,                             \eqno(28)
$$
where $Q_n=[\mu_{2n},\mu_{2n+1}]$ if $b>0,$ and
$Q_n=[\mu_{2n-1},\mu_{2n}]$ if $b<0.$
Let $\nu=\{\nu_n\}_{n\ge 0}$ be zeros of the entire
function $\theta(\pi,\la).$ Denote $\eta_n=
\mbox{sign}\,(|\theta'(\pi,\nu_n)|-|b|).$ The sequence
$\eta=\{\eta_n\}_{n\ge 0}$ is called the $\eta$- sequence
for $B.$ The inverse problem is formulated as follows

\smallskip
{\bf Inverse problem 5. } Given $\mu$ and $\eta,$ construct
$q,a$ and $b.$

\smallskip
The next theorem gives us the characterization of the spectrum
of the BVP $B.$

\smallskip
{\bf Theorem 9. }{\it For real numbers $\{\mu_n\}_{n\ge 0}\;
(\mu_n\le\mu_{n+1})$ to be the eigenvalues of a certain BVP $B$
with real potential $q(x)\in L'_2(0,\pi),$ it is necessary and
sufficient that (26) and (28) hold, where $r(\la)$ is constructed
by (27).}

\smallskip
Denote by $J'$ the set of sequences $\eta=\{\eta_n\}_{n\ge 0}$
such that\\
(i) $\eta_n=\pm 1$ if the corresponding zeros of the functions
$r(\la)$ and $r(\la)-4b$ are simple, and $\eta_n=0$ otherwise;\\
(ii) there exists $N$ (depending on the sequence) such that
$\eta_n=1$ for all $n>N.$

Clearly, if $\eta$ is the $\eta$- sequence for $B,$ then
$\eta\in J'.$ The next theorem gives us necessary and sufficient
conditions for the solvability of Inverse problem 5.

\smallskip
{\bf Theorem 10. }{\it Let real numbers $\{\mu_n\}_{n\ge 0}\;
(\mu_n\le\mu_{n+1})$ satisfying (26) and (28) be given, where
$r(\la)$ is constructed by (27). Then for each sequence $\eta\in J'$
there exists a unique real function $q(x)\in L'_2(0,\pi)$ and
real numbers $a$ and $b$ such that $\mu=\{\mu_n\}_{n\ge 0}$ is
the spectrum of $B,$ and $\eta$ is the $\eta$- sequence for $B.$}

\smallskip
We note that in [17] stability of the solution of the inverse problem
for the BVP $B$ is established.

\medskip
{\bf Acknowledgment.} This work was supported by Grant 1.1436.2014K
of the Russian Ministry of Education and Science and by Grant
13-01-00134 of Russian Foundation for Basic Research.

\begin{center}
{\bf REFERENCES}
\end{center}
\begin{enumerate}
\item[{[1]}] Marchenko V.A. Sturm-Liouville operators and their applications.
     Naukova Dumka,  Kiev, 1977;  English  transl., Birkh\"auser, 1986.
\item[{[2]}] Levitan B.M., Inverse Sturm-Liouville problems. Nauka,
     Moscow, 1984; English transl., VNU Sci. Press, Utrecht, 1987.
\item[{[3]}] P\"oschel J. and Trubowitz E., Inverse Spectral Theory,
     Academic Press, New York, 1987.
\item[{[4]}] Freiling G. and Yurko V.A., Inverse Sturm-Liouville
     Problems and their Applications. NOVA Science Publishers, New York, 2001.
\item[{[5]}] Yurko V.A. Method of Spectral Mappings in the Inverse Problem
     Theory. Inverse and Ill-posed Problems Series. VSP, Utrecht, 2002.
\item[{[6]}] Stankevich I.V., An inverse problem of spectral analysis for
     Hill's equation. Doklady Akad. Nauk SSSR 192, no.1 (1970), 34-37
     (in Russian); English transl. in Soviet Math. Dokl. 11 (1970), 582-586.
\item[{[7]}] Marchenko V.A. and Ostrovskii I.V., A characterization of the
     spectrum of the Hill operator, Mathem. Sb. 97 (1975), 540-606 (in Russian);
     English transl. in Math. USSR-Sb. 26 (1975), 4, 493-554.
\item[{[8]}] Yurko V.A. An inverse problem for second order differential operators
     with regular boundary conditions. Matem. Zametki, 18, no.4 (1975), 569-576
     (in Russian); English transl. in Mathematical Notes, 18, no.3-4 (1975), 928-932.
\item[{[9]}] Yurko V.A. On a periodic boundary value problem. Differ. Equations and
     Theory of Functions, Saratov Uni., Saratov, 1981, 109-115 (in Russian).
\item[{[10]}] Yurko V.A. On recovering differential operators with nonseparated
     boundary conditions. Study in Math. and Appl., Bashkir Uni., Ufa, 1981, 55-58
     (in Russian).
\item[{[11]}] Plaksina O.A., Inverse problems of spectral analysis for the
     Sturm-Liouville operators with nonseparated boundary conditions, Mathem. Sb. 131
     (1986), 3-26 (in Russian); English transl. in Math. USSR-Sb. 59 (1988), 1, 1-23.
\item[{[12]}] Guseinov I.M., Gasymov M.G., Nabiev I.M. An inverse problem for
     the Sturm-Liouville operator with nonseparable self-adjoint boundary conditions.
     Sib. Mathem. Zh. 31, no.6 (1990), 46-54 (in Russian); English trans. in Siberian
     Math. J. 31, no.6 (1990) 910-918.
\item[{[13]}] Guseinov I.M., Nabiev I.M. Solution of a class of inverse
     boundary-value Sturm-Liouville problems. Mathem. Sb. 186 (1995), no.5, 35-48
     (in Russian); English transl. in Sbornik: Mathematics 186 (1995), no.5, 661-674.
\item[{[14]}] Kargaev P.; Korotyaev E. The inverse problem for the Hill operator,
     a direct approach. Invent. Math. 129, no.3 (1997), 567—593.
\item[{[15]}] Yurko V.A. On differential operators with nonseparated boundary
     conditions. Funkt. Analiz i Prilozh. 28, no.4 (1994), 90-92 (Russian); English
     transl. in Functional Analysis and Applications 28, no.4 (1994), 295-297.
\item[{[16]}] Yurko V.A. The inverse  spectral  problem for differential operators with
     nonseparated boundary conditions. Journal of Mathematical Analysis and Applications
     250, no.1 (2000), 266-289.
\item[{[17]}] Freiling G., Yurko V.A. On the stability of constructing a potential in
     the central symmetry case. Applicable Analysis, 90, no.12 (2011), 1819-1828.
\end{enumerate}

\end{document}